\documentclass[showpacs,showkeys,floats,aps,nofootinbib,preprint]{revtex4}
\usepackage{bm,amsfonts, mathtools}
\usepackage{graphicx,color, wasysym}
\usepackage{textcomp}
\usepackage{amsmath,amssymb,latexsym,epsfig}
\def\nn{\nonumber}

\begin{document}

\makeatletter

\title{Chebyshev polynomials and generalized complex numbers}

\author{D. Babusci}
\email{danilo.babusci@lnf.infn.it}
\affiliation{INFN - Laboratori Nazionali di Frascati, via E. Fermi, 40, IT 00044 Frascati (Roma), Italy}

\author{G. Dattoli}
\email{giuseppe.dattoli@enea.it}

\author{E. Di Palma}
\email{emanuele.dipalma@enea.it}

\author{E. Sabia}
\email{elio.sabia@enea.it}
\affiliation{ENEA - Centro Ricerche Frascati, via E. Fermi, 45, IT 00044 Frascati (Roma), Italy}

\begin{abstract}
The generalized complex numbers can be realized in terms of $2\times2$ or higher-order matrices and can be exploited to get different 
ways of looking at the trigonometric functions. Since Chebyshev polynomials are linked to the power of matrices and to trigonometric 
functions, we take the quite natural step to discuss  them in the context of the theory of generalized complex numbers. We also briefly 
discuss the two-variable Chebyshev polynomials and their link with the third-order Hermite polynomials. 
\end{abstract}

\maketitle

We formulate the theory of Chebyshev polynomials \cite{Andrews} as a by-product of the generalized complex numbers \cite{Babusci}. 
To this aim we remind that any number $h$ satisfying the identity
\begin{equation}
\label{compl}
h^2 =  a_2 + b_2\,h
\end{equation}
is a generalized complex number\footnote{$h$ reduces to the ordinary imaginary unit for $a_2 = -1$, $b_2 = 0$.}, whose higher-order 
powers are given by 
\begin{equation}
\label{hord}
h^n = a_n + h\,b_n \qquad\qquad (n \in \mathbb{Z}^*) 
\end{equation}
with $b_n$, $a_n$ obtained recursively ($a_0 = b_1 = 1$, $a_1 = b_0 = 0$). By multiplying both sides of this equation by $h$ and equating 
the ``real" and ``imaginary" parts, i.e.,  the coefficients of the terms independent from and proportional to $h$, respectively, we obtain 
($a_2 = a$, $b_2 = b$)
\begin{equation}
\label{abrec}
a_{n + 1} = a\,b_n, \qquad\qquad b_{n + 1} = a_n + b\,b_n,
\end{equation}
which can also be rewritten in the matrix form
\begin{equation}
\label{qmat}
\left(
\begin{array}{c}
   a_{n + 1}    \\
   b_{n + 1}    
\end{array}
\right) = \hat{Q} (a, b)\,\left(
\begin{array}{c}
   a_n    \\
   b_n    
\end{array}
\right), \qquad\qquad \hat{Q} (a, b) = \left(
\begin{array}{cc}
  0  &  a  \\
  1  &  b 
\end{array}
\right),
\end{equation}
where the matrix $\hat{Q}$ satisfies identity \eqref{hord}.

From Eq. \eqref{hord} it also follows that 
\begin{equation}
e^{h\,\phi} = \sum_{n = 0}^\infty \frac{\phi^n}{n!}\,h^n = C (\phi) + h\,S (\phi),
\end{equation}
with
\begin{equation}
\label{cslik}
C (\phi) = \sum_{n = 0}^\infty \frac{\phi^n}{n!}\,a_n, \qquad\qquad
S (\phi) = \sum_{n = 0}^\infty \frac{\phi^n}{n!}\,b_n,  
\end{equation}
that provides an Euler-like identity from which, taking into account Eq. \eqref{abrec}, the following differential equations for the cos- and sin-like 
functions $C$ and $S$, respectively, can be derived
\begin{equation}
\frac{\mathrm{d}}{\mathrm{d} \phi} C (\phi) = a\,S (\phi), \qquad\qquad 
\frac{\mathrm{d}}{\mathrm{d} \phi} S (\phi) = C (\phi) + b\,S (\phi).
\end{equation}

The complex unit $h$ is characterized by the two conjugated forms (the solutions of Eq. \eqref{compl})
\begin{equation}
h_\pm = \frac{b \pm \sqrt{b^2 + 4\,a}}2
\end{equation}
and, therefore,
\begin{equation}
h_\pm^n = a_n + h_\pm\,b_n,
\end{equation}
from which we obtain 
\begin{align}
\label{abn}
a_n &= \frac{h_+\,h_-^n - h_-\,h_+^n}{h_+ - h_-} = \frac{a}{\sqrt{b^2 + 4\,a}}\,(h_+^{n - 1} - h_-^{n - 1}) \nn \\
b_n &= \frac{h_+^n - h_-^n}{h_+ - h_-} =  \frac{h_+^n - h_-^n}{\sqrt{b^2 + 4\,a}}.
\end{align}
Analogously, we define the following Euler-like identity
\begin{equation}
e^{h_\pm\,\phi} = C (\phi) + h_\pm\,S (\phi)
\end{equation}
from which we get
\begin{equation}
C (\phi) = h_+\,e^{h_+\,\phi} + h_-\,e^{h_-\,\phi}, \qquad\qquad 
S (\phi) = \frac{e^{h_+\,\phi} + e^{h_-\,\phi}}{h_+ + h_-}.
\end{equation}

Let us now consider a recurrence of the type
\begin{equation}
\label{Cherec}
L_{n + 1} = a\,L_n + b\,L_{n - 1}.
\end{equation}
By using the Binet method \cite{Yamaleev}, i.e., setting $L_n = h^n$, the solution of this difference equation can be written as 
\begin{equation}
L_n = c_+\,h_+^n + c_-\,h_-^n,
\end{equation}
with $c_\pm$ constants determined by the values of $L_0$ and $L_1$. The Chebyshev polynomials satisfy the recurrence \eqref{Cherec} with 
$a = 2\,x$ and $b = -1$ \cite{Andrews}, and, therefore, to the corresponding difference equation we can associate the unit $H$ defined by the 
equation
\begin{equation}
\label{Hpm}
H^2 = 2\,x\,H - 1, \qquad\qquad H_\pm = x \pm \sqrt{x^2 - 1}
\end{equation}
and
\begin{equation}
\label{Hun}
H_\pm^n = A_n\, + H_\pm\,B_n, 
\end{equation}
where $A_n$ and $B_n$ are functions of the variable $x$ satisfying recurrences that can be obtained by applying the same method used before 
(see Eqs. \eqref{abrec}, \eqref{qmat}), which yield
\begin{equation}
\left(
\begin{array}{c}
   A_{n + 1}    \\
   B_{n + 1}    
\end{array}
\right) = \hat{Q} (-1 , 2\,x)\,\left(
\begin{array}{c}
   A_n    \\
   B_n    
\end{array}
\right).
\end{equation}
Since $H_+ + H_- = 2\,x$, $H_+\,H_- = 1$, from Eq. \eqref{Hun} it's easy to show that
\begin{align}
  \frac{H_+^n + H_-^n}2 &= A_n + x\,B_n   \nn \\
  H_+^n\,H_-^n &= A_n^2 + 2\,x\,A_n\,B_n + B_n^2 = 1.    
\end{align} 
By replacing $H_\pm$ to $h_\pm$ in Eq. \eqref{abn}, by using again Eq. \eqref{Hpm} one obtains
\begin{align}
\label{AnBn}
A_n &= \frac{H_+\,H_-^n - H_-\,H_+^n}{H_+ - H_-} = - \frac{H_+^{n - 1} - H_-^{n - 1}}{2\,\sqrt{x^2 - 1}} \nn \\
B_n &= \frac{H_+^n - H_-^n}{H_+ - H_-} = - A_{n + 1}.
\end{align}
Furthermore, since 
\begin{equation}
\label{dHpm}
\frac{\mathrm{d} H_\pm}{\mathrm{d} x} = \pm \frac{H_\pm}{\sqrt{x^2 -1}}
\end{equation}
we find that functions $B_n$ verify the following differential equation
\begin{equation}
\left[(1 - x^2)\,\partial_x^2 - 3\,x\,\partial_x + (n - 1)^2\right] B_n (x) = 0
\end{equation}
from which, by performing the replacement  $n \to n+1$, we obtain the differential equation for the Chebyshev polynomials of second kind $U_n (x)$, 
and we can therefore make the following identification\footnote{The Chebyshev polynomials of first kind are given by 
$T_n (x) =  (H_+^n + H_-^n)/2$.} 
\begin{equation}
U_n (x) =  \frac{H_+^{n + 1} - H_-^{n + 1}}{2\,i\,\sqrt{1 - x^2}}.
\end{equation}
By specifying Eq. \eqref{hord} to the case of matrix $\hat{Q}$ for the Chebyshev polynomials, taking into account Eqs. \eqref{AnBn}, \eqref{Cherec}, 
one has
\begin{equation}
\hat{Q}^{n + 1} (1, - 2\,x) = A_{n + 1} + \hat{Q} (1, -2\,x)\,B_{n + 1} = \left(
\begin{array}{cc}
  - U_{n - 1} (x)  & - U_n (x)  \\
    U_n (x)  &  U_{n + 1} (x)   
\end{array}
\right),
\end{equation}
from which, since $\det \hat{Q} (1, - 2\,x) = 1$, we get the well-known identity
\begin{equation}
U_n^2 (x) - U_{n - 1} (x)\,U_{n + 1} (x) = 1.
\end{equation}

It is well known that any $2\times2$ matrix can be written as a linear combination of the Pauli matrices and the unit matrix as follows \cite{BDD}
\begin{equation}
\hat{M} = \alpha\,\hat{1}_2 + \sum_{k = 1}^3 \beta_k\,\hat{\sigma}_k.
\end{equation}
As a consequence of the identity $\left\{\hat{\sigma}_i, \hat{\sigma}_j\right\} = 2\,\delta{jk}$, we get 
\begin{equation}
\hat{M}^2 = \gamma\,\hat{1}_2 + 2\,\alpha\,\hat{M} \qquad\qquad (\gamma = - \alpha^2 + \sum_{k = 1}^3 \beta_k^2),
\end{equation}
and, therefore, according to previous point of view, any matrix $\hat{M}$ can be viewed as the realization of a generalized complex unit. Furthermore, 
it can also be checked that \cite{Ricci1, Dattoli}  
\begin{equation}
\hat{M}^n = U_{n - 1} (\alpha)\,\hat{M} + U_{n - 2} (\alpha)\,\hat{1}_2.
\end{equation}
A step forward in the theory of Chebyshev polynomials has been done in \cite{Ricci2} where their two-variable generalization $U^{(2)}_n (u, v)$ has 
been introduced. These polynomials satisfy the recurrence 
\begin{equation}
U^{(2)}_{n + 2} = u\,U^{(2)}_{n + 1} - v\,U^{(2)}_n + U^{(2)}_{n - 1} 
\end{equation}
and, therefore, are associated with the (third-order) imaginary unit defined by
\begin{equation}
Y^3 = u\,Y^2 - v\,Y + 1
\end{equation}
and to higher-order trigonometry discussed in Refs. \cite{Babusci} and \cite{Yamaleev}. 

We close the paper by noting that the generating function of the two-variable Chebyshev polynomials of second kind is given by \cite{Ricci2}
\begin{equation}
\sum_{n = 0}^\infty t^n\,U^{(2)}_{n + 1} (u, v) = \frac1{1 - u\,t + v\,t^2 - t^3}.
\end{equation}
The use of Laplace transform yields
\begin{equation}
\frac1{1 - u\,t + v\,t^2 - t^3} = \int_0^\infty \mathrm{d}s\,e^{- s\,(1 - u\,t + v\,t^2 - t^3)}
\end{equation}
and, therefore, on account of the generating function of the third-order Hermite polynomials
\begin{equation}
\sum_{n = 0}^\infty \frac{t^n}{n!}\,H_n^{(3)} (x, y, z) = e^{x\,t + y\,t^2 + z\,t^3},
\end{equation}
we get the following identity
\begin{equation}
U^{(2)}_{n + 1} (u, v) = \frac1{n!}\,\int_0^\infty \mathrm{d}s\,e^{- s}\,H_n^{(3)} (u\,s, - v\,s, s),
\end{equation}
that can significantly simplifies the analysis of the properties of the two-variable Chebyshev polynomials, and of their multivariable generalization 
as well. 


\begin{thebibliography}{99}

\bibitem{Andrews} L. C. Andrews, \textit{Special Functions for Engineers and Applied Mathematicians}, Mc Millan, New York (1985).

\bibitem{Babusci} D. Babusci, G. Dattoli, E. Di Palma, and E. Sabia, Adv. Appl. Clifford Al. {\bf 22}, 271 (2012).

\bibitem{Seroul} R. S\'eroul, \textit{Programming for Mathematicians}, p. 21, Springer-Verlag, Berlin (2000). 

\bibitem{Yamaleev} R. M. Yamaleev, Adv. Appl. Clifford Al. {\bf 10}, 15 (2005).

\bibitem{BDD} D. Babusci, G. Dattoli, and M. Del Franco, \textit{Lectures on Mathematical Methods for Physics}, 
Internal Report ENEA RT/2010/5837 (available at www.frascati.enea.it/biblioteca/).

\bibitem{Ricci1} P. E. Ricci, Atti Accad. Sc. Torino, {\bf 109} (1974-75), and references therein.

\bibitem{Dattoli} G. Dattoli, C. Mari, and A. Torre, Nuovo Cimento B {\bf 108}, 61 (1993).

\bibitem{Ricci2} P. E. Ricci, Rendiconti di Matematica, Serie VI, vol. {\bf 11}, 295 (1978).

\end{thebibliography}
\end{document}